\documentclass[a4paper,12pt]{article}
\usepackage{amsmath,amsfonts,amssymb,amscd}

\long\def\bee#1\ee{\begin{equation}#1\end{equation}}
\def\be{\begin{equation}}

\def\bea{\begin{eqnarray}}
\def\eea{\end{eqnarray}}

\def\s{\scriptstyle }
\def\a{\alpha}
\def\b{\beta}

\def\L{\Lambda}

\def\N{{\mathbb N}}

\def\A{{\mathbb A}}

\def\cE{{\mathcal E}}

\def\cD{{\mathcal D}}

\def\moins{\raise 1pt\hbox{$\s -$}}
\def\plus{\raise 1pt \hbox{ $\s +$}}

\newtheorem{theorem}{Theorem}

\newtheorem{lemma}[theorem]{Lemma}

\begin{document}

\begin{center} 
\bf\Large The differential equation \\
satisfied by a plane curve of degree $\mathbf n$
\end{center}

\medskip\centerline{Alain Lascoux}

\bigskip
\begin{abstract}
Eliminating the arbitrary coefficients in the equation of 
a generic plane curve of order $n$ by computing sufficiently many
derivatives, one obtains a differential equation.
This is a projective invariant. The first one,
corresponding to conics, has been obtained by Monge. 
Sylvester, Halphen, Cartan used invariants of higher order.
The expression of these invariants is rather complicated, 
but becomes much simpler when interpreted in terms
of symmetric functions.

\medskip
\centerline{\bf R\'esum\'e}

\smallskip
L'expression diff\'erentielle des courbes planes de degr\'e donn\'e
fournit un invariant projectif. Monge a obtenu celle des coniques 
planes, Sylvester et Halphen ont g\'en\'eralis\'e l'\'equation de
Monge aux courbes planes de tout degr\'e. Nous montrons que la
th\'eorie des fonctions sym\'etriques permet de retrouver ces 
invariants, et d'en donner des expressions plus compactes.

\end{abstract}
\medskip

\noindent{\bf Mots-clefs}: \'Equation de Monge, invariant d'Halphen,
invariants diff\'erentiels.


\bigskip
 
A line in the plane can be written
$$ y = a\, x +b \ ,$$
with arbitrary coefficients $a,b$, but it is more satisfactory 
to write it
$$   y'' =0   \ .$$

More generally 
$$
  { y}^n + ({ \star}) { y}^{n-1} x+\cdots +
 ({ \star}) x^n +\cdots + ({ \star}) { y}^0 x^0 = 0 $$
is the equation of a general planar curve of order $n$.
Writing sufficiently many derivatives of this equation, one can 
eliminate in their system the arbitrary coefficients $({ \star})$.

However, already in the case of a conic (solved by Monge),
we have to use the derivatives of order $3,4,5$, and the outcome
is not straightforward to interpret. We need some method
to perform the elimination.

It is convenient, instead of taking $\frac{d^i}{dx^i}$,  to 
rather use the normalization
$$   D^i  = \frac{d^i}{i!\, dx^i}$$

With these conventions, the Leibnitz formula loses its coefficients:
$$ D^n(fg) = \sum_{i+j=n}  D^if\, D^j g \ .$$

We need the collection $\{ D^0y,\, D^1y,\, D^2y, \ldots \}$, 
which we can write with the help of a generating series
$$  \sum_{i=0}^\infty   z^i\,  D^iy \ .$$

Symmetric  function theory tells that we ought to formally
factorize this series, as we factorize the total Chern class of a vector
bundle. 

Thus we introduce a formal  alphabet $\A$ and write 
\begin{equation}
\sum_{i=0}^\infty \frac{d^i y}{i!\, dx^i}  
    = \sum_{i=0}^\infty   z^i\,  D^iy = 
\prod_{a\in \A} (1+za) = \sum_{i=0}^\infty z^i\, \L^i\A \ ,
\end{equation}
denoting by $\L^i\A$ the \emph{elementary symmetric functions} in 
$\A$, and thus identifying $D^iy$ to $\L^i\A$. We refer to
\cite{Mac} for what concerns the theory of symmetric functions,
and to \cite{Cbms} for its $\lambda$-ring approach.  

Remember that taking $k$ copies of an alphabet 
(we write $k\,\A$) translates into 
taking the $k$-th power of the generating function~: 
$$  \left(\sum_{i=0}^\infty z^i\, \L^i\A   \right)^k  =
 \sum_{i=0}^\infty z^i\, \L^i(k\, \A)  \ .$$
Adding $r$ copies of a letter $x$ 
to these $k$ copies of $\A$ is written, at the level
of generating series, as
\begin{equation}  \label{eq:fg}
  (1+zx)^r\, \left(\sum_{i=0}^\infty z^i\, \L^i\A   \right)^k  =
 \sum_{i=0}^\infty z^i\, \L^i(k\, \A +r\, x)  \ .
\end{equation}

Thus, instead of having a sum  
$ D^n(y^2) = \sum_{i+j=n}  D^iy\, D^j y$ to express the derivatives 
of the square of $y$, 
one can now  use the more compact notation
$D^n(y^2) = \L^n(2\A)$. 

More generally,
 $D^n(y^3) = \L^n(3\A)$, $D^n(y^4) = \L^n(4\A),\ldots $ 
and one has the following 
easy lemma resulting from Leibnitz' rule:

\begin{lemma}
Given $n,k,r\in\N$, then  
\begin{equation}  \label{eq:der}
  D^n( x^r y^k)  =  x^r \L^n ( k\A + r/x)   \ .
\end{equation}
\end{lemma}

We can now easily write the derivatives of any orders of the 
components of the equation of a planar curve.

Let us look first at the case treated by Monge. 

We start with 
$$ u =  y^2 + c_1 xy +c_2 y + (\star) x^2 + (\star) x +(\star)  $$
and take successive derivatives, starting from the third one
(so that the part depending on $x$ only (where the coefficients
$(\star)$ appear)  has already been eliminated).
\begin{eqnarray} 
 D^3 u &=&  \L^3(2\A) + x c_1 \L^3 (\A+1/x) +  c_2 \L^3\A \\
       &=& \L^3(2\A) +(xc_1+c_2) \L^3(\A)   +  c_1 \L^2\A \\
 D^4 u &=& \L^4(2\A) +(xc_1+c_2) \L^4(\A)   +  c_1 \L^3\A \\ 
 D^5 u &=& \L^5(2\A) +(xc_1+c_2) \L^5(\A)   +  c_1 \L^4\A 
\end{eqnarray}
Elimination of the coefficients among  these three equations gives
the vanishing~:
\begin{equation}
  \begin{vmatrix}   \L^2(\A) & \L^3(\A) & \L^3(2\A) \\
                    \L^3(\A) & \L^4(\A) & \L^4(2\A) \\
                    \L^4(\A) & \L^5(\A) & \L^5(2\A) \\
   \end{vmatrix}
 = \begin{vmatrix}   \L^2 & \L^3 & 2\L^{30} + 2\L^{21} \\
                    \L^3 & \L^4 & 2\L^{40} +2\L^{31}+\L^{22} \\
                    \L^4 & \L^5 &  2\L^{50} +2\L^{41}+2\L^{32} 
   \end{vmatrix}  \ ,
\end{equation}
writing $\L^i$ for $\L^i\A$, $\L^{ij}$ 
for $\L^i\A \L^j\A,\ldots$ .

The last determinant can be simplified and becomes
\begin{equation}   \label{eq:Conic2}
  \begin{vmatrix}   \L^2 & \L^3 & 0 \\
                    \L^3 & \L^4 & \L^{22} \\
                    \L^4 & \L^5 & 2\L^{32}
   \end{vmatrix}  =
\begin{vmatrix}   y''/2 & y'''/6 & 0 \\
                    y'''/6 & y^{iv}/24 & (y''/2)^2 \\
                    y^{iv}/24 & y^{v}/120 & y'' y'''/6
   \end{vmatrix}   \ ,
\end{equation}
which is Monge's equation, after suppressing the extra factor $y''/2$~: 
\begin{equation}
 D^2y D^2y D^5y - 3 D^2y D^3y D^4y + 2 D^3y D^3y D^3y  = 0  \ .
\end{equation}

The general case takes only a few lines  more to be written down.

From Eq.\ref{eq:fg}, one has 
$$\L^n(\A+ rx) = \L^n\A + rx\, \L^{n-1}\A + \binom{r}{2}x^2\, \L^{n-2}\A +
          \cdots +  \binom{r}{n}x^n\, \L^{0}\A   \ .$$
Let $u$ be a polynomial in $x,y$ of total degree $n$ with leading
term  $y^n$.  The equations 
$$ D^{n+1}u= 0 = \cdots = D^{n(n+3)/2} u $$
are
\begin{eqnarray*}
0 &=& \L^{n+1}(n\A) +c_{1,n-1} \L^{n}((n\moins 1)\A)
                      + c_{2,n-1} \L^{n+1}((n\moins1)\A) +  \cdots\\ 
 & &\hspace{140pt}
         +\,  c_{1,1} \L^2\A +c_{2,1} \L^3\A+\cdots + \L^{n+1}\A \\
0 &=& \L^{n+2}(n\A) +c_{1,n-1} \L^{n+1}((n\moins 1)\A)+ c_{2,n-1} \L^{n+2}((n\moins 1)\A)
   +  \cdots\\ & &\hspace{140pt}
            +\, c_{1,1} \L^3\A +c_{2,1} \L^4\A+\cdots + \L^{n+2}\A \\
  & & \cdots   \qquad \cdots \qquad  \cdots \qquad \cdots\qquad \cdots\qquad \cdots  \\
0 &=& \L^{n(n+3)/2}(n\A) \\ 
  & & \hspace{40pt}   +\, c_{1,n-1} \L^{n(n+3)/2-1}((n\moins 1)\A)
   + c_{2,n-1} \L^{n(n+3)/2}((n\moins 1)\A)+  \cdots    \\
 & & \hspace{60pt} +\,   c_{1,1} \L^{n(n+1)/2+1}\A 
      +c_{2,1} \L^{n(n+1)/2+2}\A+\cdots + \L^{n(n+3)/2}\A \ , 
\end{eqnarray*}
where the coefficients $c_{i,j}$ are polynomials in $x$ only.

Eliminating these coefficients, one obtains the vanishing of the
following determinant (we have written the columns in a different order):
\begin{equation}    \label{eq:PremiereEq}
 \begin{vmatrix}
  \L^2\A  & \L^3\A & \cdots & \L^{n+1}\A &\cdots & \L^{n}((n\moins 1)\A)&
    \L^{n+1}((n\moins 1)\A)  & \L^{n+1}(n\A)  \\
 \L^3\A  & \L^4\A & \cdots & \L^{n+2}\A &\cdots & \L^{n+1}((n\moins 1)\A)&
    \L^{n+2}((n\moins 1)\A)  & \L^{n+2}(n\A)  \\
  \vdots & \vdots &     &\vdots &   & \vdots & \vdots &\vdots  \\
   & & \cdots & \L^N\A & \cdots & \L^{N-1}((n\moins 1)\A)
  & \L^N((n\moins 1)\A) & \L^N(n\A)
 \end{vmatrix}
\end{equation} 
with $N= n(n+3)/2$, 
which is the differential equation satisfied by a planar curve of order $n$.

This determinant has a simple structure, with blocks of $n,\, n\moins1,\
\ldots,\, 1$ columns involving the elementary symmetric functions 
in $\A,\,  2\A,\, \ldots,\, n\A$.

One can simplify it a little. Because the image of a curve of degree
$n$ under the transformation $y \to \a y+\b$ 
is still a curve of the same degree,
the value of the determinant is independent of $y=\L^0\A$ and 
$y'=\L^1\A$, that one can put both equal to $0$. 

Therefore, instead of the generating series (\ref{eq:fg}), 
one can now take a second alphabet $\cD$ such that
$$ \L^i\cD= \L^{i+2}\A \ ,\ i=0,1,\ldots $$
(as usual $\L^i=0$ for $i<0$).

In other words, 
$\sum_i z^i\, \L^i\cD =  \L^2\A +z \L^3\A +\cdots$, and for what
concerns its powers, one has that
$$  \L^i(k\cD)=  \L^{i+2k}(k\A) \quad , \quad i,k\in\N \ .  $$

The equation of a planar curve can now be rewritten
\begin{equation}    \label{eq:DeuxiemeEq}
\hspace{-10pt} \begin{vmatrix}
  \L^0\cD  &  \cdots & \L^{n-1}\cD  & \L^{-1}(2\cD)&\cdots  & \L^{n-3}(2\cD)
               & \cdots & \L^{1-n}(n\cD) \\
  \L^1\cD  &  \cdots & \L^{n}\cD  & \L^{0}(2\cD)&\cdots &  \L^{n-2}(2\cD)
               & \cdots & \L^{2-n}(n\cD) \\
  \vdots  &   &\vdots &\vdots & &\vdots  &   & \vdots \\
   &\cdots &  \L^N\cD  & &\cdots & \L^{N-2}(2\cD)   & &
     \L^{N-2n+2}(n\cD)
 \end{vmatrix}   
\end{equation}
with $N= (n-1)(n+4)/2$, 
and, apart from notations, this is the equation given by Sylvester.

For the conic, this equation is the determinant (\ref{eq:Conic2})
that we have seen above~:
\begin{equation}     \begin{vmatrix} 
   \L^0\cD  & \L^1\cD  & \L^{-1}(2\cD)  \\
   \L^1\cD  & \L^2\cD  & \L^{0}(2\cD)  \\
    \L^2\cD  & \L^3\cD  & \L^{1}(2\cD)  
\end{vmatrix}  =
 \begin{vmatrix}
   \L^0\cD  & \L^1\cD  & 0  \\
   \L^1\cD  & \L^2\cD  & \L^{0}\cD\, \L^{0}\cD   \\
    \L^2\cD  & \L^3\cD  & 2\L^{1}\cD\,\L^{0}\cD
\end{vmatrix}
\end{equation}

The equation of a planar cubic is~:
\begin{equation}   \label{eq:cubic}    
 \begin{vmatrix}
\L^0\cD &\L^1\cD  & \L^2\cD    & \L^{-1}(2\cD) & \L^{0}(2\cD) & \L^{-2}(3\cD)   \\
 \L^1\cD &\L^2\cD  & \L^3\cD   & \L^{0}(2\cD) & \L^{1}(2\cD) & \L^{-1}(3\cD)   \\
 \vdots & \vdots   & \vdots    &\vdots       &\vdots   &\vdots  \\
\L^5\cD &\L^6\cD  & \L^7\cD    & \L^{4}(2\cD) & \L^{5}(2\cD) & \L^{3}(3\cD)   
\end{vmatrix} = 0 \ .  
\end{equation}
In terms of $\L^i\cD$ only, written $\L_i$, and putting $\L_0=1$,
 the determinant reads
$$
\begin{vmatrix}  1& \Lambda_{1}& \Lambda_{2}& 0& 1& 0\cr \Lambda_{
1}& \Lambda_{2}& \Lambda_{3}& 1& 2\Lambda_{1}& 0\cr \Lambda_{
2}& \Lambda_{3}& \Lambda_{4}& 2\Lambda_{1}& \Lambda_{1}^{2}+
2\Lambda_{2}& 1\cr \Lambda_{3}& \Lambda_{4}& \Lambda_{5}&
\Lambda_{1}^{2}+2\Lambda_{2}& 2\Lambda_{1}\Lambda_{2}+2\Lambda_{
3}& 3\Lambda_{1}\cr \Lambda_{4}& \Lambda_{5}& \Lambda_{6}
& 2\Lambda_{1}\Lambda_{2}+2\Lambda_{3}& \Lambda_{2}^{2}+2\Lambda_{
1}\Lambda_{3}+2\Lambda_{4}& 3\Lambda_{1}^{2}+3\Lambda_{2}\cr
\Lambda_{5}& \Lambda_{6}& \Lambda_{7}& \Lambda_{2}^{2}+2\Lambda_{
1}\Lambda_{3}+2\Lambda_{4}& 2\Lambda_{2}\Lambda_{3}+2\Lambda_{
1}\Lambda_{4}+2\Lambda_{5}& \Lambda_{1}^{3}+6\Lambda_{1}\Lambda_{
2}+3\Lambda_{3}
\end{vmatrix} $$
which expands into the rather less attractive expression~:
$6\L_{1}^{7}\L_{2}\L_{3}^{2}-30\L_{2}^{4}
\L_{1}^{2}\L_{5}-10\L_{3}^{3}\L_{1}\L_{
5}+5\L_{2}\L_{7}\L_{3}^{2}-2\L_{5}\L_{
7}\L_{1}^{3}+5\L_{7}\L_{1}^{2}\L_{2}^{3}
+30\L_{4}\L_{2}^{5}\L_{1}+6\L_{5}\L_{
2}^{2}\L_{3}^{2}+6\L_{1}^{5}\L_{6}\L_{4}
-6\L_{1}^{6}\L_{6}\L_{3}+3\L_{6}\L_{
1}\L_{4}^{2}-3\L_{6}^{2}\L_{1}\L_{2}-20\L_{
1}^{4}\L_{4}^{2}\L_{3}-60\L_{4}\L_{1}^{3
}\L_{2}^{4}+5\L_{2}\L_{5}\L_{4}^{2}+10\L_{
5}\L_{4}\L_{3}^{2}-3\L_{1}^{8}\L_{5}\L_{
2}+20\L_{2}^{2}\L_{4}^{2}\L_{3}-3\L_{4}\L_{
1}\L_{5}^{2}+10\L_{2}^{2}\L_{4}^{2}\L_{1
}^{3}-15\L_{7}\L_{1}^{4}\L_{2}^{2}+2\L_{
6}^{2}\L_{1}^{3}+4\L_{7}\L_{2}^{4}+5\L_{
2}^{6}\L_{3}-8\L_{2}\L_{4}\L_{3}\L_{
1}\L_{5}-5\L_{2}\L_{6}\L_{1}^{2}\L_{
5}-10\L_{3}^{2}\L_{2}^{2}\L_{1}\L_{4}+5\L_{
1}^{3}\L_{4}^{3}+12\L_{1}^{7}\L_{5}\L_{3
}-6\L_{2}\L_{1}^{5}\L_{4}^{2}+5\L_{1}^{2
}\L_{4}\L_{3}^{3}+5\L_{1}^{2}\L_{3}\L_{
5}^{2}-30\L_{1}^{4}\L_{2}\L_{3}^{3}-50\L_{
2}^{4}\L_{3}^{2}\L_{1}-12\L_{5}\L_{1}^{6
}\L_{2}^{2}+6\L_{2}^{2}\L_{4}\L_{1}^{7}-
10\L_{4}\L_{2}^{4}\L_{3}-10\L_{4}\L_{
2}^{3}\L_{5}-10\L_{1}^{6}\L_{3}^{3}+10\L_{
2}\L_{4}\L_{1}^{4}\L_{5}-7\L_{1}^{2}\L_{
6}\L_{4}\L_{3}+12\L_{2}\L_{7}\L_{1}
^{6}+10\L_{2}^{3}\L_{3}^{3}+\L_{1}^{2}\L_{
5}\L_{4}^{2}+3\L_{3}\L_{6}\L_{1}\L_{
5}-3\L_{3}\L_{7}\L_{1}\L_{4}-2\L_{4
}\L_{6}\L_{5}-15\L_{2}\L_{6}\L_{4}\L_{
1}^{3}+5\L_{2}\L_{7}\L_{1}^{2}\L_{4}-2\L_{
2}\L_{6}\L_{4}\L_{3}+10\L_{3}^{2}\L_{
1}\L_{4}^{2}+6\L_{1}^{7}\L_{4}^{2}-10\L_{
2}\L_{1}\L_{4}^{3}-10\L_{1}\L_{2}\L_{
3}^{4}+3\L_{5}\L_{7}\L_{1}\L_{2}-5\L_{
1}^{3}\L_{6}\L_{3}^{2}+10\L_{1}^{4}\L_{5
}\L_{3}^{2}-3\L_{1}^{8}\L_{4}\L_{3}+30\L_{
1}^{5}\L_{4}\L_{3}^{2}-6\L_{1}^{6}\L_{4}
\L_{5}-4\L_{1}^{5}\L_{7}\L_{3}-2\L_{
1}^{5}\L_{5}^{2}-10\L_{6}\L_{2}^{4}\L_{1
}-6\L_{6}\L_{2}^{3}\L_{3}-4\L_{7}\L_{
2}^{2}\L_{4}+4\L_{6}\L_{2}^{2}\L_{5}+30\L_{
6}\L_{2}^{3}\L_{1}^{3}-3\L_{1}^{8}\L_{7}
+54\L_{2}^{2}\L_{3}^{2}\L_{1}^{5}+50\L_{
2}^{2}\L_{3}^{3}\L_{1}^{2}+24\L_{4}\L_{1
}^{5}\L_{2}^{3}-5\L_{6}\L_{3}^{3}-30\L_{
2}^{3}\L_{1}^{3}\L_{3}^{2}+45\L_{2}^{3}\L_{
1}^{4}\L_{5}+60\L_{2}^{5}\L_{1}^{2}\L_{3
}-45\L_{2}^{4}\L_{1}^{4}\L_{3}-10\L_{2}^{
3}\L_{1}\L_{4}^{2}-10\L_{2}^{3}\L_{1}^{6
}\L_{3}-5\L_{3}\L_{4}^{3}+3\L_{2}^{5}\L_{
1}^{5}+10\L_{2}\L_{1}^{2}\L_{4}^{2}\L_{3
}-14\L_{7}\L_{1}\L_{2}^{2}\L_{3}+10\L_{
2}^{6}\L_{1}^{3}+5\L_{2}\L_{1}^{3}\L_{5}
^{2}-27\L_{6}\L_{1}^{5}\L_{2}^{2}-15\L_{
2}^{7}\L_{1}+40\L_{2}^{3}\L_{3}\L_{1}\L_{
5}+6\L_{2}\L_{6}\L_{1}^{7}-8\L_{2}\L_{
3}\L_{5}^{2}+6\L_{3}^{5}+\L_{5}^{3}+\L_{
1}^{9}\L_{6}+\L_{6}^{2}\L_{3}-20\L_{2}\L_{
4}\L_{3}^{3}+10\L_{4}\L_{1}^{2}\L_{2}^{3
}\L_{3}+45\L_{4}\L_{1}^{4}\L_{3}\L_{
2}^{2}-40\L_{5}\L_{1}^{3}\L_{3}\L_{2}^{2
}-15\L_{6}\L_{1}^{2}\L_{3}\L_{2}^{2}+14\L_{
6}\L_{4}\L_{1}\L_{2}^{2}+10\L_{2}\L_{
7}\L_{3}\L_{1}^{3}+18\L_{2}\L_{6}\L_{
1}\L_{3}^{2}-5\L_{2}\L_{5}\L_{1}^{2}\L_{
3}^{2}-42\L_{2}\L_{4}\L_{1}^{6}\L_{3}-20
\L_{2}\L_{4}\L_{1}^{3}\L_{3}^{2}-12\L_{
2}\L_{5}\L_{1}^{5}\L_{3}+20\L_{2}\L_{
6}\L_{1}^{4}\L_{3}+\L_{7}\L_{4}^{2}+\L_{
1}^{2}\L_{7}\L_{3}^{2}-\L_{5}\L_{7}\L_{
3}$
Recall that 
$\L_0= y''/2!$, $\L_1= y'''/6!,\ldots,\, \L_7= y^{vii}/7!$, so that 
the equation of a cubic is a differential polynomial of order $15$, and degree $10$. 

\bigskip
\centerline{\bf \Large Invariance under the projective group}

\bigskip 
We did not yet use that 
the image of a curve under a projective transformation of the plane
is still a curve of the same degree, and therefore that 
\emph{the differential equation is a projective invariant}.

Expressed in terms of $\cD$, invariance under some subgroup of the projective
group amounts to belonging  to the kernel of
$$\nabla_{\cD}:=  \frac{d}{d_{\L^1}} + 2\L^1\cD \frac{d}{d_{\L^2}} +
    3\L^2\cD \frac{d}{d_{\L^3}}+  4\L^3\cD \frac{d}{d_{\L^4}} +\cdots  $$

This operator appears in the theory of \emph{binary forms} in $x,y$, and
simply expresses  the invariance under the translation $x\to x+1$.
Elements of this kernel are called \emph{semi-invariants} in the theory
of binary forms.

To simplify the operator $\nabla_{\cD}$, let
us introduce another alphabet $\cE$~:
\begin{equation}
 \L^i\cE=  \L^i\cD/i!  = \frac{1}{i!\, (i\plus 2)!}
                                \frac{d^{i+2} y}{dx^{i+2}}  \ .
\end{equation}

Under the change of alphabet, our projective invariants
belong to the kernel of
$$\nabla_{\cE}:=  \frac{d}{d_{\L^1}} + \L^1\cE \frac{d}{d_{\L^2}} +
    \L^2\cE \frac{d}{d_{\L^3}}+  \L^3\cE \frac{d}{d_{\L^4}} +\cdots  $$

But this kernel is very easy to determine. Indeed, 
 $\nabla_{\cE}$ sends $\psi_1\cE$ onto $1$, 
and it sends all the other power sums $\psi_i\cE$, $i=2,\ldots,n$  
 onto $0$, $n$ being the cardinality of $\cE$.
Therefore a semi-invariant is a polynomial in 
$\psi_2\cE,\,  \psi_3\cE, \ldots $.

 Monge's element is of degree $3$ in $\cE$. Therefore it must be
proportional to $\psi_3\cE$.  Indeed, the equation of a conic 
is
\begin{equation}
   \fbox{$\,  \psi_3\, \cE  = 0 \, $}  \ , 
\end{equation}
or, in terms of elementary symmetric functions, 
\begin{equation}
   \psi_3\cE = 3\Lambda^{3}\cE-3\Lambda^{12}\cE +\Lambda^{111}\cE   
  = 0 
\end{equation}
which rewrites, with $\L^1\cE= \frac{d^{3} y}{1!\, 3!!\,  dx^{3}}$,
$\L^2\cE= \frac{d^{4} y}{2!\, 4!!\,  dx^{4}}$, 
$\L^3\cE= \frac{d^{5} y}{1!\, 3!!\,  dx^{5}}$, 
into the equation that we have already seen:
$$ { 3   \left(\frac{y''}{2!}\right)^2
\frac{y^{v}}{3!5!} - 3
\frac{y''}{2!} \frac{y'''}{1!3!} \frac{y^{iv}}{2!4!}
+ \left(\frac{y'''}{1!3!}\right)^3 = 0 \ .} $$

The equations of curves of higher degree can be similarly treated,
using the usual theory of binary forms, and involve
generalizations of the Hessian.

The next projective invariant, after the invariant of Monge,
has been found by Halphen\cite{Halphen}. 
Cartan\cite{Cartan} takes it as the projective analogue of the curvature.
The equation of the cubic is a polynomial in it and Monge
invariant. Halphen invariant, in terms of power sums in $\cE$, is

\begin{equation}
48\psi_{5}\psi_{3}-20\psi_{3}^{2}\psi_{2}-\psi_{2}^{4}+12\psi_{
2}^{2}\psi_{4}-36\psi_{4}^{2}
\end{equation}
and is easier to handle than determinants of the type
displayed in Eq. (\ref{eq:cubic}).

\vspace{2cm}
\begin{center}
Alain Lascoux \\
CNRS, Institut Gaspard Monge, Universit\'e de Marne-la-Vall\'ee\\
77454 Marne-la-Vall\'ee Cedex, France\\[2pt]
{\tt Alain.Lascoux@univ-mlv.fr}
\end{center}

\end{document}